\documentclass[10pt,reqno]{amsart}
\usepackage{hannah}
\usepackage{color,graphicx}
\usepackage[latin1]{inputenc}

\newcommand {\brk}{\rule {1mm}{0mm}}
\newcommand {\df}[1]{\emph {#1}}

\newcommand {\calP}{{\mathcal P}}

\DeclareMathOperator{\trop}{trop}
\DeclareMathOperator{\cplx}{cplx}

\DeclareMathOperator{\Area}{Area}
\DeclareMathOperator{\path}{path}
\DeclareMathOperator{\irr}{irr}

\title [The Caporaso-Harris formula in tropical geometry]{The Caporaso-Harris
  formula and plane relative Gromov-Witten invariants in tropical geometry}
\author{Andreas Gathmann and Hannah Markwig}
\address {Andreas Gathmann, University of Kaiserslautern, Department of
  Mathematics, PO Box 3049, 67653 Kaiserslautern, Germany}
\email {andreas@mathematik.uni-kl.de}
\address {Hannah Markwig, University of Kaiserslautern, Department of
  Mathematics, PO Box 3049, 67653 Kaiserslautern, Germany}
\email {markwig@mathematik.uni-kl.de}
\thanks {The second author has been funded by the DFG grant Ga 636/2.}

\begin{document}

\begin{abstract}
  Some years ago Caporaso and Harris have found a nice way to compute the
  numbers $ N(d,g) $ of complex plane curves of degree $d$ and genus $g$
  through $ 3d+g-1 $ general points with the help of relative Gromov-Witten
  invariants. Recently, Mikhalkin has found a way to reinterpret the
  numbers $ N(d,g) $ in terms of tropical geometry and to compute them by
  counting certain lattice paths in integral polytopes. We relate these two
  results by defining an analogue of the relative Gromov-Witten invariants and
  rederiving the Caporaso-Harris formula in terms of both tropical geometry and
  lattice paths.
\end{abstract}

\maketitle


\section{Introduction}

Let $ N(d,g) $ be the number of complex curves of degree $d$ and genus $g$ in
the projective plane $ \P^2 $ through $ 3d+g-1 $ general fixed points. These
numbers have first been computed by Caporaso and Harris \cite{CH98}. Their
strategy is to define ``relative Gromov-Witten invariants'' that count plane
curves of given degree and genus having specified local contact orders to a
fixed line $L$ and passing in addition through the appropriate number of
general points. By specializing one point after the other to lie on $L$ they
derive recursive relations among these relative Gromov-Witten invariants that
finally suffice to compute all the numbers $ N(d,g) $.

As an example of what happens in this specialization process we consider plane
rational cubics having a point of contact order 3 to $L$ at a fixed point $ p_1
\in L $ and passing in addition through 5 general points $ p_2,\dots,p_6 \in
\P^2 $. To compute the number of such curves we move $ p_2 $ to $L$. What
happens to the cubics under this specialization? As they intersect $L$ already
with multiplicity 3 at $ p_1 $ they cannot pass through another point on $L$
unless they become reducible and have $L$ as a component. There are two ways
how this can happen: they can degenerate into a union of three lines $ L \cup
L_1 \cup L_2 $ where $ L_1 $ and $ L_2 $ each pass through two of the points $
p_3,\dots,p_6 $, or they can degenerate into $ L \cup C $, where $C$ is a conic
tangent to $L$ and passing through $ p_3,\dots,p_6 $:

\begin {center} \input {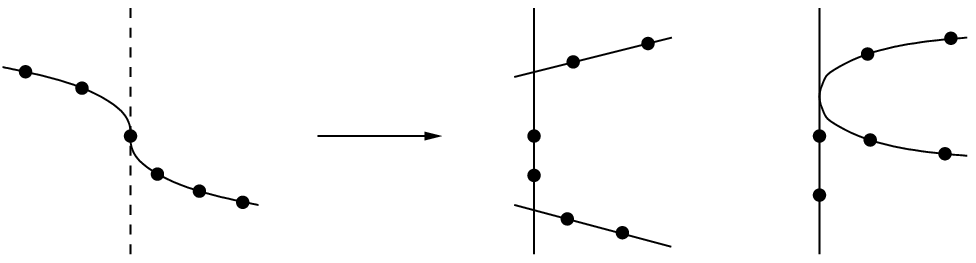} \end {center}

The initial number of rational cubics with a point of contact order 3 to $L$ at
a fixed point and passing through 5 more general points is therefore a sum of
two numbers (counted with suitable multiplicities) related to only lines and
conics. This is the general idea of Caporaso and Harris how specialization
finally reduces the degree of the curves and allows a recursive solution to
compute the numbers $ N(d,g) $.

Recently, Mikhalkin found a different way to compute the numbers $ N(d,g) $
using tropical geometry \cite{Mi03}. He proved the so-called ``Correspondence
Theorem'' that asserts that $ N(d,g) $ can be reinterpreted as the number of
\emph{tropical} plane curves of degree $d$ and genus $g$ through $ 3d+g-1 $
points in general position. The goal of this paper is to establish a connection
between the complex and the tropical point of view. We show that relative
Gromov-Witten invariants correspond to tropical curves with unbounded edges of
higher weight to the left. For example, our cubics above having a point of
contact order 3 to $L$ in a fixed point correspond to tropical curves having a
fixed unbounded end of weight 3 to the left (and passing in addition through 5
more general points $ p_2,\dots,p_6 $). The specialization process above is
simply accomplished by moving $ p_2 $ to the very far left:

\begin {center} \input {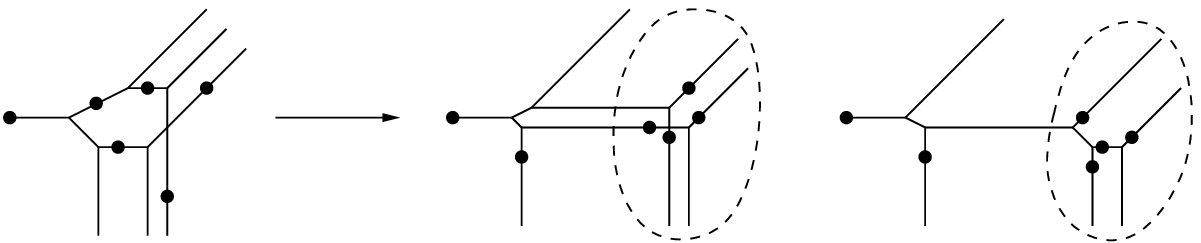} \end {center}

We see that the resulting tropical curves ``split'' into two parts: a left part
(through $ p_2 $) and a right part (through the remaining points, circled in
the picture above). We get the same ``degenerations'' as in the complex case:
one where the right part consists of two lines through two of the points $
p_3,\dots,p_6 $ each, and one where it consists of a conic ``tangent to a
line'' (i.e.\ with an unbounded edge of multiplicity 2 to the left).

To arrive at the actual numbers $ N(d,g) $ Mikhalkin did not count these
tropical curves directly however. Instead, he showed by purely combinatorial
arguments that the number $ N(d,g) $ is also equal to the number of certain
``increasing'' lattice paths of length $ 3d+g-1 $ from $ (0,d) $ to $ (d,0) $
in the integral triangle $ \Delta_d = \{ (x,y) \in \Z^2 ;\; x \ge 0, y \ge 0,
x+y \le d \} $, counted with suitable multiplicities. We will show that the
idea of Caporaso and Harris can also be seen directly in this lattice path
set-up: relative Gromov-Witten invariants simply correspond to lattice paths
with fixed integral steps on the left edge of $ \Delta_d $. For example, the
cubics above with triple contact to $L$ correspond to lattice paths starting
with the two points $ (0,3) $ and $ (0,0) $, i.e.\ with a step of length 3. The
remaining steps are then arbitrary, as long as the number of steps in the path
is correct (in this case 6 as we have 6 marked points):

\begin {center} \input {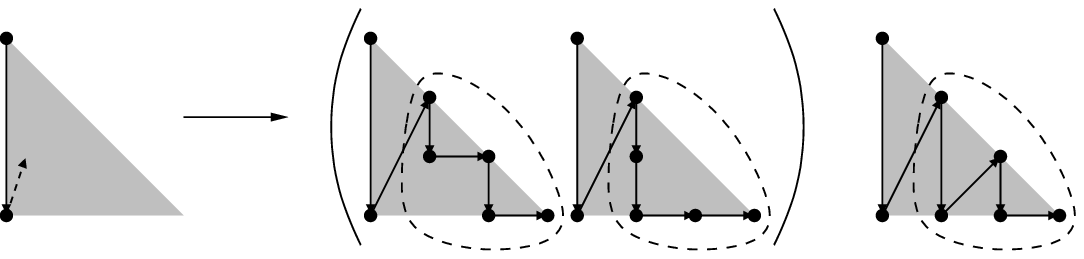} \end {center}

If we look at the triangle $ \Delta_{d-1} $ obtained from $ \Delta_d $ by
removing the left edge (circled in the picture above) we see again two possible
types: one that corresponds to a union of two lines (the first two cases)
and one that corresponds to conics tangent to $L$ (the last case, with a
step of length 2 at the left side of $ \Delta_{d-1} $).

The aim of this paper is to make the above ideas rigorous. We will define
``relative Gromov-Witten invariants'' both in terms of tropical curves and
lattice paths, and prove the Caporaso-Harris formula in both settings. For
simplicity we will work with not necessarily irreducible curves most of the
time (except for section \ref{ch-irr}). We will also restrict to complex curves
in the plane (i.e.\ tropical curves with Newton polyhedron $ \Delta_d $) to
keep the notation simple. However, the same ideas can be applied for other
toric surfaces as well (see e.g.\ remarks \ref {rem-gen} and \ref {ch-rect}).
We hope that our ideas can also be generalized to curves in higher-dimensional
spaces. Work on this question is in progress.

This paper is organized as follows. In section \ref{repeat} we will brief\/ly
review the known results on complex and tropical curves. We will then construct
analogues of relative Gromov-Witten invariants and prove the Caporaso-Harris
formula in the lattice path set-up in section \ref{sec-lattice}. The same is
then done for the tropical curves set-up in section \ref{sec-tropical}.

We would like to thank Ilia Itenberg for pointing out a serious error in a
previous version of this paper.


\section{Complex and tropical curves} \label{repeat}

In this section we will brief\/ly review the notations and results on complex
and tropical curves that we need in our paper. Our main references are
\cite{CH98} and \cite{Mi03}.

\subsection{Complex curves and the Caporaso-Harris formula}

We start by defining the ``relative Gromov-Witten invariants'' used by Caporaso
and Harris to compute the numbers of complex plane curves of given degree and
genus through the appropriate number of given points.

\begin {definition}
  A \df {(finite) sequence} is a collection $ \alpha=(\alpha_1,\alpha_2,\dots)
  $ of natural numbers almost all of which are zero. If $ \alpha_k=0 $ for all
  $ k>n $ we will also write this sequence as $
  \alpha=(\alpha_1,\dots,\alpha_n) $. For two sequences $ \alpha $ and $ \beta
  $ we define
    \[ \begin {array}{r@{\;}l@{\;}l}
      |\alpha| &:=& \alpha_1+\alpha_2+\cdots; \\
      I\alpha &:=& 1\alpha_1+2\alpha_2+3\alpha_3+\cdots; \\
      I^{\alpha} &:=& 1^{\alpha_1}\cdot 2^{\alpha_2}\cdot 3^{\alpha_3}\cdot
        \; \cdots ; \\
      \alpha+\beta &:=& (\alpha_1+\beta_1,\alpha_2+\beta_2,\ldots); \\
      \alpha \geq \beta &:\Leftrightarrow& \alpha_n \geq \beta_n
        \mbox { for all $n$}; \\
      \binom{\alpha}{\beta} &:=& \binom{\alpha_1}{\beta_1}\cdot
        \binom{\alpha_2}{\beta_2}\cdot \; \cdots.
    \end {array} \]
  We denote by $e_k$ the sequence which has a $1$ at the $k$-th place and zeros
  everywhere else.
\end {definition}

\begin {definition} \label{def-cplx}
  Let $ d \ge 0 $ and $g$ be integers, and let $ \alpha $ and $ \beta $ be two
  sequences with $ I\alpha+I\beta=d $. Pick a fixed line $ L \subset \P^2 $.
  Then we denote by $ N_{\cplx}^{\alpha,\beta}(d,g) $ the number of smooth but
  not necessarily irreducible curves (or more precisely: stable maps) of degree
  $d$ and genus $g$ to $ \P^2 $ that
  \begin {itemize}
  \item intersect $L$ in $ \alpha_i $ fixed general points of $L$ with
    multiplicity $i$ each for all $ i \ge 1 $;
  \item intersect $L$ in $ \beta_i $ more arbitrary points of $L$ with
    multiplicity $i$ each for all $ i \ge 1 $; and
  \item pass in addition through $ 2d+g+|\beta|-1 $ more general points in $
    \P^2 $.
  \end {itemize}
  In other words, we consider the numbers of complex plane curves of given
  degree and genus that have fixed contact orders to a given line. Note that
  this set of numbers includes the numbers $ N_{\cplx} (d,g) :=
  N_{\cplx}^{(0),(d)} (d,g) $ of complex plane curves of degree $d$ and genus
  $g$ through $ 3d+g-1 $ general points (without multiplicity conditions).
\end {definition}

The main result of Caporaso and Harris is how these numbers can be computed
recursively:

\begin {definition}
  We say that a collection of numbers $ N^{\alpha,\beta} (d,g) $ defined for
  all integers $ d \ge 0 $ and $g$ and all sequences $ \alpha,\beta $ with $
  I\alpha+I\beta =d $ \df{satisfies the Caporaso-Harris formula} if
  \begin{align*}
    N^{\alpha,\beta}(d,g)
      =& \sum_{k:\beta_k>0} k\cdot N^{\alpha+e_k,\beta-e_k}(d,g)
        \\
      &\qquad + \sum I^{\beta'-\beta} \cdot \binom{\alpha}{\alpha'} \cdot
        \binom{\beta'}{\beta} \cdot N^{\alpha',\beta'}(d-1,g')
  \end{align*}
  for all $ d,g,\alpha,\beta $ as above with $ d>1 $, where the second sum is
  taken over all $\alpha'$, $\beta'$ and $g'$ satisfying 
  \begin{align*}
    \alpha' &\leq \alpha; \\
    \beta' &\geq \beta; \\
    I\alpha'+I\beta'&=d-1; \\
    g-g'&=|\beta'-\beta|-1; \\
    d-2 &\geq g-g'.
  \end{align*}
\end {definition}

\begin {theorem} \label{thm-ch}
  The numbers $ N_{\cplx}^{\alpha,\beta} (d,g) $ of definition \ref{def-cplx}
  satisfy the Caporaso-Harris formula.
\end {theorem}

\begin {proof}
  See \cite{CH98}.
\end {proof}

Note that a collection of numbers $ N^{\alpha,\beta}(d,g) $ satisfying the
Caporaso-Harris formula is determined uniquely by their values for $ d=1 $. In
particular, theorem \ref{thm-ch} allows us to compute all numbers $
N_{\cplx}(d,g) $ from the starting information that there is exactly one line
through two points in the plane.

\subsection{Tropical curves}

We will move on to the tropical set-up and recall the definition of tropical
curves and their Newton polyhedra (following \cite{Mi03} and \cite {NS04}).

Let $\overline{\Gamma}$ be a weighted finite graph without divalent vertices.
The weight of an edge $E$ of $\overline{\Gamma}$ will be written as $\omega(E)
\in \N \setminus \{0\}$. Denote the set of $1$-valent vertices by
$\overline{\Gamma}^0_{\infty}$. We remove the one-valent vertices from
$\overline{\Gamma}$ and set $\Gamma:= \overline{\Gamma}\setminus
\overline{\Gamma}^0_{\infty}$. The graph $\Gamma$ may then have non-compact
edges which are called \df{unbounded edges} or \df{ends}. We write $\Gamma^0$
for the set of vertices and $\Gamma^1$ for the set of edges of $ \Gamma $.  We
define the set of \df{flags} of $\Gamma$ by $\F \Gamma := \{ (V,E) \in \Gamma^0
\times \Gamma^1 \mid V \in \partial E\}$. Let $ h:\Gamma \rightarrow \R^2 $ be
a continuous proper map such that the image $h(E)$ of any edge $E \in \Gamma^1$
is contained in an affine line with rational slope. Then we can define a map
$u:\F \Gamma \rightarrow \Z^2$ that sends $(V,E)$ to the primitive integer
vector that starts at $h(V)$ and points in the direction of $h(E)$.

\begin{definition}
  A \df{parametrized tropical (plane) curve} is a pair $ (\Gamma,h) $ as above
  such that
  \begin{enumerate}
  \item for every edge $E \in \Gamma^1$ the restriction $h\vert_E$ is an
    embedding;
  \item for every vertex $V \in \Gamma^0$ the \df{balancing condition}
    \begin{displaymath}
      \sum_{E \in \Gamma^1: V \in \partial E} \omega (E) \cdot u(V,E) =0.
    \end{displaymath} 
    holds.
  \end{enumerate}
  A \df{tropical curve} in $ \R^2 $ is the image $ h(\Gamma) $ of a
  parametrized tropical curve.
\end{definition}

Let $V$ be an $r$-valent vertex of a tropical curve $\Gamma$ and let
$E_1,\ldots, E_r$ be the counterclockwise enumerated edges adjacent to $V$.
Draw in the $\Z^2$-lattice an orthogonal line $l(E_i)$ of integer length
$\omega(E_i)$ to $h(E_i)$, where $l(E_1)$ starts at any lattice point and
$l(E_i)$ starts at the endpoint of $l(E_{i-1})$, and where by ``integer
length'' we mean $\sharp (\Z^2 \cap l(E_i)) -1$. The balancing condition tells
us that we end up with a closed $r$-gon. If we do this for every vertex we end
up with a polygon in $\Z^2$ that is divided into smaller polygons. The polygon
is called the \df{Newton polygon} of the tropical curve, and the division the
corresponding \df{Newton subdivision}. Note that the ends of the curve
correspond to lines on the boundary of the Newton polygon.

\begin {example}
  The following picture shows an example of a tropical curve and its Newton
  polygon and subdivision:

  \begin {center} \begin{picture}(0,0)%
\includegraphics{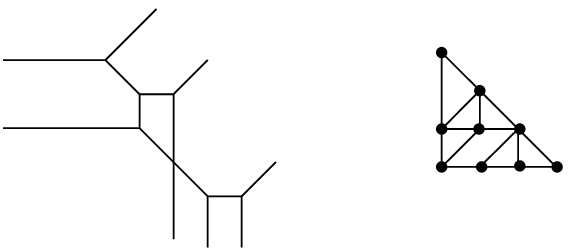}%
\end{picture}%
\setlength{\unitlength}{3947sp}%
\begingroup\makeatletter\ifx\SetFigFont\undefined%
\gdef\SetFigFont#1#2#3#4#5{%
  \reset@font\fontsize{#1}{#2pt}%
  \fontfamily{#3}\fontseries{#4}\fontshape{#5}%
  \selectfont}%
\fi\endgroup%
\begin{picture}(2705,1168)(3661,-992)
\put(3676,-40){\makebox(0,0)[lb]{\smash{{\SetFigFont{8}{9.6}{\familydefault}{\mddefault}{\updefault}{\color[rgb]{0,0,0}$\omega=2$}%
}}}}
\end{picture}%
 \end {center}
\end {example}

Most properties of tropical curves can be read off both from their image in $
\R^2 $ as well as in the dual picture from their Newton polygon and
subdivision. Here are some examples:

\begin{definition} \label{tropdef} \brk
  \vspace {-1.4\baselineskip}
  \begin{enumerate}
  \item A parametrized tropical plane curve is called \df{simple} if $\Gamma$
    is $3$-valent, the map $h$ is injective on the set of vertices, for a
    vertex $V$ and an edge $E$ we have $ h(V)\cap h(E)=\emptyset $, for two
    edges $E_1$ and $E_2$ we have $\#\{h(E_1)\cap h(E_2)\}\leq 1 $, and for all
    $ p \in \R^2 $ we have $ \sharp h^{-1}(p)\leq 2 $. A tropical plane curve
    is called simple if it allows a simple parametrization. In the dual
    language, a curve is simple if and only if its subdivision contains only
    triangles and parallelograms.
  \item The \df{genus} of a parametrized tropical curve $ (\Gamma,h) $ is
    defined to be the number $ 1-\dim H^0(\Gamma,\Z)+ \dim H^1(\Gamma,\Z) =
    1-\#\Gamma^0+\#\Gamma^1_0 $. The genus of a simple tropical plane
    curve is the least genus of all parametrizations that the curve allows. In
    the dual of a simple tropical curve, the genus is the number of interior
    lattice points of the subdivision minus the number of parallelograms.
  \item A parametrized tropical curve is called \emph {irreducible} if the
    graph $\Gamma$ is connected. A simple tropical plane curve is called
    \emph{irreducible} if it allows only irreducible parametrizations. An
    \emph{irreducible component} of a simple tropical plane curve is the image
    of a connected component of the graph $\Gamma$ of a parametrization with
    the maximum possible number of connected components.
  \item The \df{degree} $\Delta$ of a tropical plane curve is the collection of
    directions of the unbounded edges together with the sum of the weights for
    each direction. In the dual language the degree is just the Newton polygon,
    which is therefore also denoted by $\Delta$.
  \item The \df{(combinatorial) type} of a parametrized tropical curve is given
    by the weighted graph $\Gamma$ together with the map $u:\F \Gamma
    \rightarrow \Z^2$. The (combinatorial) type of a tropical curve is the
    combinatorial type of any parametrization of least possible genus. In the
    dual setting, the type is just the Newton subdivision.
  \item \label{tropdef-e}
    Let $V$ be a trivalent vertex of $\Gamma$ and $E_1,E_2,E_3$ the edges
    adjacent to $V$. The \df{multiplicity} of $V$ is the product of the area of
    the parallelogram spanned by $u(V,E_1)$ and $u(V,E_2)$ times the weights of
    the edges $E_1$ and $E_2$. The balancing condition tells us that this
    definition is independent of the order of $E_1,E_2$ and $E_3$. In the dual
    language, the multiplicity of $3$-valent vertex is equal to $2$ times the
    area of the dual triangle.
  \item \label{tropdef-f}
    The \df{multiplicity} $\mult(C)$ of a simple tropical plane curve is the
    product of the multiplicities of all trivalent vertices of $\Gamma$ as in
    \ref{tropdef-e}. In the dual language, the multiplicity is the product over
    all double areas of triangles of the dual subdivision.
  \end{enumerate}
\end{definition}

\begin {example}
  The following picture shows an irreducible simple tropical curve (on the
  left) and a simple tropical curve with two components (on the right),
  together with their parametrizations of least possible genus. Both curves are
  of degree $ (3 \cdot (0,-1), 3 \cdot (-1,0),3 \cdot (1,1)) $.

  \begin {center} \input {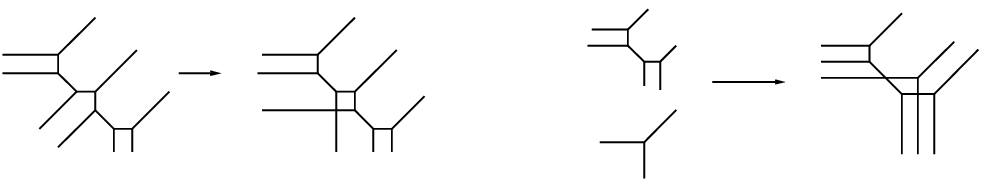} \end {center}
\end {example}

We are now ready to state the results of \cite{Mi03} that are relevant for our
purposes.

\begin {definition}
  Let $ d \ge 0 $ and $g$ be integers, and let $ \Delta_d $ be the Newton
  polyhedron $ \{ (x,y) \in \Z^2 ;\; x\ge 0, y \ge 0, x+y\le d \} $. Let $
  \calP = \{p_1,\dots,p_{3d+g-1}\} \subset \R^2 $ be a set of points in
  (tropical) general position (see \cite{Mi03} section 4.2 for a precise
  definition). Then by \cite{Mi03} lemma 4.22 the number of tropical curves
  through all points of $ \calP $ is finite, and all such curves are simple.
  We denote by $ N_{\trop} (d,g) $ the number of tropical curves of degree $
  \Delta_d $ and genus $g$ passing through all points of $ \calP $, counted
  with their multiplicities as in definition \ref{tropdef} \ref{tropdef-f}. By
  \cite{GM051} this definition does not depend on the choice of points.
\end {definition}

\begin{remark} \label{rem-strings}
  Let $\mathcal{P}=(p_1,\ldots,p_n)$ be a configuration of points in general
  position, and let $(\Gamma,h)$ be a (simple) tropical curve through
  $\mathcal{P}$. Then by \cite{Mi03} lemma 4.20 each connected component of
  $ \Gamma \setminus \{h^{-1}(p_1),\ldots, h^{-1}(p_n)\} $ is a tree and
  contains exactly one unbounded edge. That is, in a tropical curve through
  $\mathcal{P}$, we can neither find a path which connects two unbounded edges
  nor a path around a loop without meeting a marked point. 
\end{remark}

\begin {theorem}[``Correspondence Theorem''] \label{corrthm}
  For all $d$ and $g$ we have $ N_{\trop}(d,g) = N_{\cplx}(d,g) $.
\end {theorem}

\begin {proof}
  See \cite {Mi03} theorem 1.
\end {proof}

\subsection{Lattice paths} \label{lpath}

In his paper \cite{Mi03} Mikhalkin also gave an algorithm to compute the
numbers $ N_{\trop}(d,g) $ combinatorially. He did not calculate them directly
however, but rather by a trick that relates them to certain lattice paths that
we will now introduce.

\begin{definition}
  A path $\gamma: [0,n] \rightarrow \R^2$ is called a lattice path if
  $\gamma |_{[j-1,j]}, j=1,\ldots,n$ is an affine-linear map and $\gamma(j)
  \in\Z^2 $ for all $ j=0\ldots,n $.
\end{definition}

Let $ \lambda $ be a fixed linear map of the form $ \lambda: \R^2 \to \R,
\lambda (x,y) = x-\varepsilon y $, where $ \varepsilon $ is a small irrational
number. A lattice path is called $\lambda$-increasing if $\lambda \circ \gamma$
is strictly increasing. Let $p:=(0,d) $ and $q:=(d,0) $ be the points in $
\Delta:=\Delta_d $ where $\lambda|_{\Delta}$ reaches its minimum (resp.\
maximum). The points $p$ and $q$ divide the boundary $\partial \Delta$ into two
$\lambda$-increasing lattice paths $\delta_{+}:[0,n_+]\rightarrow \partial
\Delta$ (going clockwise around $\partial \Delta$) and $ \delta_-: [0,n_-]
\rightarrow \partial \Delta$ (going counterclockwise around $\partial \Delta$),
where $n_{\pm}$ denotes the number of integer points in the $\pm$-part of the
boundary. The following picture shows an example for $ d=3 $:

\begin {center} \input {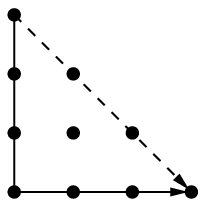} \end {center}

We will now define the multiplicity of $ \lambda $-increasing paths as in
\cite{Mi03}:

\begin{definition} \label{def-mu}
  Let $\gamma:[0,n]\rightarrow \Delta$ be a $\lambda$-increasing path from $p$
  to $q$, that is, $\gamma(0)=p$ and $\gamma(n)=q$. The \df{multiplicities}
  $\mu_+(\gamma)$ and $\mu_-(\gamma)$ are defined recursively as follows:
  \begin{enumerate}
  \item \label{def-mu-a}
    $\mu_{\pm}(\delta_{\pm}):=1$.
  \item \label{def-mu-b}
    If $\gamma \neq \delta_{\pm}$ let $k_{\pm} \in [0,n]$ be the smallest
    number such that $\gamma$ makes a left turn (respectively a right turn for
    $\mu_-$) at $\gamma(k_{\pm})$. (If no such $k_{\pm}$ exists we set
    $\mu_{\pm}(\gamma):=0$). Define two other $\lambda$-increasing lattice
    paths as follows:
    \begin {itemize}
    \item $\gamma_{\pm}':[0,n-1]\rightarrow \Delta$ is the path that cuts the
      corner of $\gamma(k_{\pm})$, i.e.\ $\gamma'_{\pm}(j):=\gamma(j)$ for
      $j<k_{\pm}$ and $\gamma'_{\pm}(j) := \gamma(j+1)$ for $j \geq k_{\pm}$.
    \item $\gamma''_{\pm}:[0,n]\rightarrow \Delta$ is the path that completes
      the corner of $\gamma(k_{\pm})$ to a parallelogram, i.e.\ $\gamma''_{\pm}
      (j):=\gamma(j)$ for all $j \neq k_{\pm}$ and $\gamma''_{\pm}(k_{\pm})
      :=\gamma(k_{\pm}-1)+\gamma(k_{\pm}+1)-\gamma(k_{\pm})$:

      \begin {center} \input {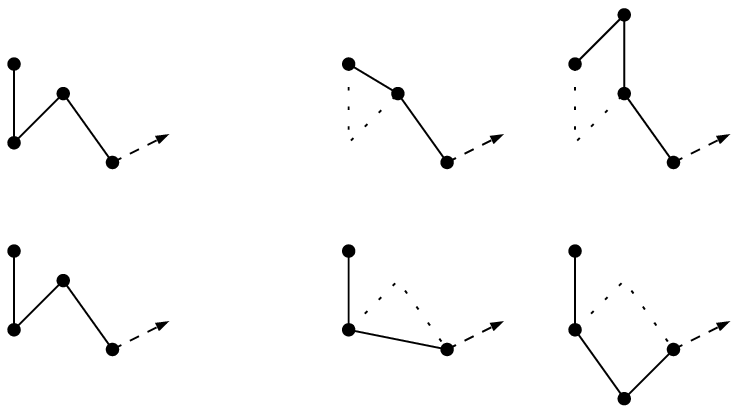} \end {center}

    \end {itemize}
    Let $T$ be the triangle with vertices $\gamma(k_{\pm}-1),\gamma(k_{\pm}),
    \gamma(k_{\pm}+1)$. Then we set
    \begin{equation*}
      \mu_{\pm}(\gamma):=2 \cdot \Area{T} \cdot \mu_{\pm}(\gamma'_{\pm})
        +\mu_{\pm}(\gamma''_{\pm}).
    \end{equation*}
    As both paths $\gamma'_{\pm}$ and $\gamma''_{\pm}$ include a smaller area
    with $\delta_{\pm}$, we can assume that their multiplicity is known. If
    $\gamma''_{\pm}$ does not map to $\Delta$, $\mu_{\pm}(\gamma''_{\pm})$ is
    defined to be zero.
  \end{enumerate}
  Finally, the multiplicity $\mu(\gamma)$ is defined to be the product
  $\mu(\gamma) := \mu_+(\gamma) \mu_-(\gamma)$.
\end{definition}

Note that the multiplicity of a path $\gamma$ is positive only if the recursion
above ends with the path $\delta_+:[0,n_+]\rightarrow \Delta$ (respectively
$\delta_-$). In other words, if we end up with a ``faster'' path
$\delta':[0,n']\rightarrow \Delta$ such that
$\delta_+([0,n_+])=\delta'([0,n'])$ but $n'<n_+$ then the multiplicity is zero.

\begin{definition}
  Let $ d \ge 0 $ and $g$ be integers. We denote by $ N_{\path} (d,g) $ the
  number of $ \lambda $-increasing lattice paths $\gamma:[0,3d+g-1]
  \rightarrow \Delta$ with $\gamma(0)=p$ and $\gamma(3d+g-1)=q$ counted with
  their multiplicities as in definition \ref {def-mu}.
\end{definition}

\begin{theorem} \label{corrthm2}
  For all $d$ and $g$ we have $N_{\path}(d,g)=N_{\trop}(d,g)$.
\end{theorem}

\begin {proof}
  See \cite{Mi03} theorem 2. The idea is to choose a line $H$ orthogonal to the
  kernel of $\lambda$ and $n=3d+g-1$ points $p_1,\ldots,p_n$ on $H$ such
  that the distance between $p_i$ and $p_{i+1}$ is much bigger than the
  distance of $p_{i-1}$ and $p_i$ for all $i$. These points are then in
  tropical general position. Consider a tropical curve of degree $ \Delta_d $
  and genus $g$ through these points. If we take the edges on which the marked
  points lie and consider their dual edges in the Newton subdivision then these
  dual edges can be shown to form a $ \lambda $-increasing path from $p$ to
  $q$. Furthermore, the multiplicity of a path coincides with the number of
  possible Newton subdivisions times their multiplicity. (For a path $\gamma$
  there can be several possible Newton subdivisions.) The multiplicity $\mu_+$
  counts the possible Newton subdivisions times their multiplicity in the
  half-plane above $H$, whereas $\mu_-$ counts below $H$. Passing from $\gamma$
  to $\gamma'$ and $\gamma''$ corresponds to moving the line up (for $\mu_+$)
  respectively down. The path $\gamma'$ leaves a $3$-valent vertex of
  multiplicity $2\Area(T)$ out, and $\gamma''$ counts the possibility that
  there might be a crossing of two edges, dual to a parallelogram.
  
  As an example, the following picture shows a path $\gamma$ and the possible
  Newton subdivisions for $\gamma$, all dual to tropical curves of multiplicity
  $1$. The multiplicity of the path is $3$.

  \begin {center} \begin{picture}(0,0)%
\includegraphics{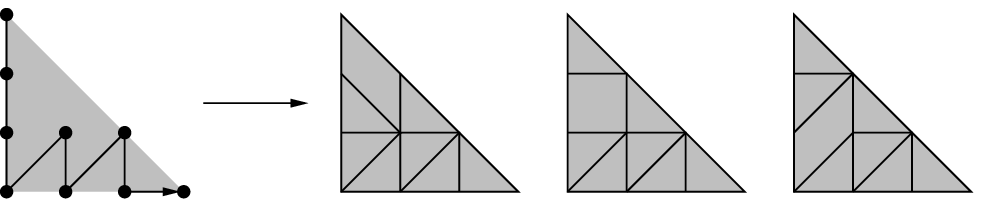}%
\end{picture}%
\setlength{\unitlength}{4144sp}%
\begingroup\makeatletter\ifx\SetFigFont\undefined%
\gdef\SetFigFont#1#2#3#4#5{%
  \reset@font\fontsize{#1}{#2pt}%
  \fontfamily{#3}\fontseries{#4}\fontshape{#5}%
  \selectfont}%
\fi\endgroup%
\begin{picture}(4452,1074)(871,-611)
\put(1306,-556){\makebox(0,0)[b]{\smash{{\SetFigFont{8}{9.6}{\familydefault}{\mddefault}{\updefault}{\color[rgb]{0,0,0}$\gamma$}%
}}}}
\end{picture}%
 \end {center}
\end {proof}

For integers $ d\ge 0 $ and $g$ we have now defined three sets of numbers $
N_{\cplx}(d,g) $, $ N_{\trop}(d,g) $, and $ N_{\path}(d,g) $. Moreover, we
know that these three sets of numbers coincide by theorems \ref{corrthm} and
\ref{corrthm2}. However, only in the complex picture have we defined numbers $
N_{\cplx}^{\alpha,\beta}(d,g) $ corresponding to curves with specified contact
orders that satisfied the Caporaso-Harris formula. It is the goal of this paper
to establish the same thing for the tropical and lattice path set-ups.


\section {The Caporaso-Harris formula in the lattice path set-up}
  \label{sec-lattice}

\subsection {Generalized lattice paths} \label {gen-lat}

We will now slightly generalize the definitions of section \ref{lpath} in order
to allow more lattice paths and arrive at lattice path analogues of the numbers
$ N_{\cplx}^{\alpha,\beta}(d,g) $. Although this concept can be used for some
other polygons (see \ref{rem-gen} and \ref{ch-rect}) we will present it here
only for the triangle $\Delta_d$ with vertices $(0,0),(0,d), (d,0)$. Choose two
sequences $ \alpha $ and $ \beta $ with $I\alpha+I\beta=d$. As above let
$\lambda(x,y)=x-\varepsilon y$ where $\varepsilon$ is a small irrational
number.

Let $\gamma:[0,n]\rightarrow \Delta_d$ be a $ \lambda $-increasing path with $
\gamma(0) = (0,d-I\alpha)=(0,I\beta)$ and $\gamma(n)=q=(d,0)$. We are going to
define a multiplicity for such a path $ \gamma $. Again this multiplicity will
be the product of a ``positive'' and a ``negative'' multiplicity that we define
separately.

\begin {definition} \label {neg-mult}
  Let $\delta_{\beta}:[0,|\beta|+d]\rightarrow \Delta_d$ be a path such that
  the image $\delta_{\beta}([0,|\beta|+d])$ is equal to the piece of boundary
  of $\Delta_d$ between $(0,I\beta)$ and $q=(d,0)$, and such that there are
  $\beta_i$ steps (i.e. images of a size one interval $[j,j+1]$) of integer
  length $i$ at the side $s$ (and hence at $\{y=0\}$ only steps of integer
  length $1$). We define the negative multiplicity $\mu_{\beta,-}
  (\delta_{\beta})$ of all such paths to be $1$. For example, the following
  picture shows all paths $ \delta_\beta $ for $ \beta=(2,1) $ and $ d=5 $:

  \begin {center} \begin{picture}(0,0)%
\includegraphics{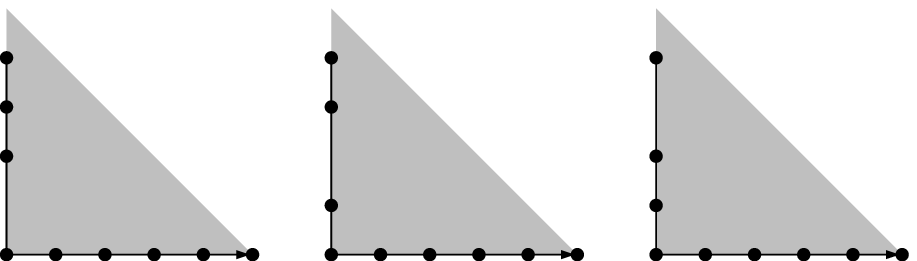}%
\end{picture}%
\setlength{\unitlength}{4144sp}%
\begingroup\makeatletter\ifx\SetFigFont\undefined%
\gdef\SetFigFont#1#2#3#4#5{%
  \reset@font\fontsize{#1}{#2pt}%
  \fontfamily{#3}\fontseries{#4}\fontshape{#5}%
  \selectfont}%
\fi\endgroup%
\begin{picture}(4155,1155)(331,-585)
\end{picture}%
 \end {center}

  Using these starting paths the \df{negative multiplicity} $ \mu_{\beta,-}
  (\gamma) $ of an arbitrary path as above is now defined recursively by the
  same procedure as in definition \ref{def-mu} \ref{def-mu-b}.
\end {definition}

\begin {definition} \label {pos-mult}
  To compute the \df{positive multiplicity} $\mu_{\alpha,+}(\gamma)$ we extend
  $\gamma$ to a path $\gamma_{\alpha}:[0,|\alpha|+n]\rightarrow \Delta_d$ by
  adding $\alpha_i$ steps of integer length $i$ at $\{x=0\}$ from
  $\gamma_{\alpha}(0)=p$ to $\gamma_{\alpha}(|\alpha|)=(0,I\beta)$. Then we
  compute $\mu_+(\gamma_{\alpha})$ as in definition \ref{def-mu} and set
  $\mu_{\alpha,+}(\gamma):= \frac{1}{I^{\alpha}}\cdot \mu_+(\gamma_{\alpha})$.
\end {definition}

\begin {remark}
  Note that definition \ref {pos-mult} seems to depend on the order in which we
  add the $ \alpha_i $ steps of lengths $i$ to the path $ \gamma $ to obtain
  the path $ \gamma_\alpha $. It will follow however from the alternative
  description of the positive multiplicity in proposition \ref {mult-formula}
  \ref {mult-formula-b} that this is indeed not the case.
\end {remark}

We can now define the analogue of relative Gromov-Witten invariants in the
lattice path set-up.

\begin{definition} \label{nd}
  Let $d \ge 0 $ and $g$ be integers, and let $ \alpha $ and $ \beta $ be
  sequences with $I\alpha+I\beta=d$. We define $N_{\path}^{\alpha,\beta}
  (d,g)$ to be the number of $\lambda$-increasing paths $\gamma:[0,2d+g+|\beta|
  -1]\rightarrow \Delta_d$ that start at $(0,d-I\alpha)=(0,I\beta)$ and end at
  $(d,0)$, where each such path is counted with multiplicity $\mu_{\alpha,
  \beta}(\gamma):=\mu_{\alpha,+}(\gamma) \cdot \mu_{\beta,-}(\gamma)$.

  Note that as expected (i.e.\ as in the complex case) we always have $
  N_{\path}(d,g) = N_{\path}^{(0),(d)} (d,g) $ by definition.
\end{definition}

\begin {example}
  The following picture shows that $ N_{\path}^{(0,1),(1)}(3,0)=4+2+1+1+2=10 $:

  \begin {center} \input {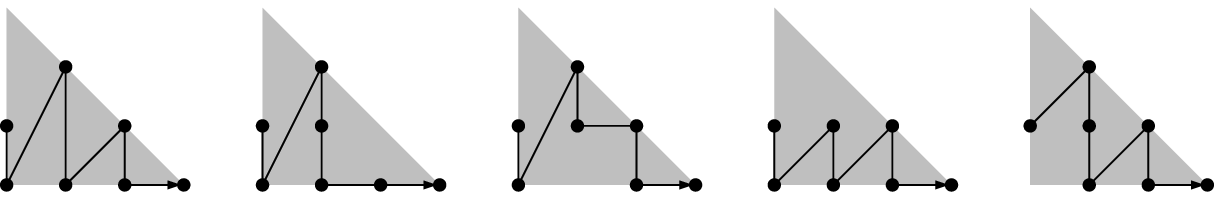} \end {center}
\end {example}

\subsection {The Caporaso-Harris formula}

In order to prove the Caporaso-Harris formula for the numbers $
N_{\path}^{\alpha,\beta} (d,g) $ of definition \ref{nd} we will first express
the negative and positive multiplicities of a generalized lattice path in a
different, non-recursive way. For this we need an easy preliminary lemma:

\begin {lemma} \label {no-skip}
  Let $ \alpha $ and $ \beta $ be two sequences with $ I\alpha+I\beta=d $, and
  let $ \gamma $ be a generalized lattice path as in section \ref {gen-lat}. If
  $ \gamma $ has a step that ``moves at least two columns to the right'', i.e.\
  that starts on a line $ \{ x=i \} $ and ends on a line $ \{ x=j \} $ for some
  $ i,j $ with $ j-i \ge 2 $ then $ \mu_{\beta,-} (\gamma) = \mu_{\alpha,+}
  (\gamma) = \mu_{\alpha,\beta} (\gamma) = 0 $.
\end {lemma}

\begin {proof}
  If $ \gamma $ is a path with a step that moves at least two columns to the
  right then the same is true for the paths $ \gamma_{\pm}' $ and $
  \gamma_{\pm}'' $ of definition \ref {def-mu}. Hence the lemma follows by
  induction since the only end paths $ \delta_\beta $ (see definition \ref
  {neg-mult}) and $ \delta_+ $ (see definition \ref {def-mu}) with non-zero
  multiplicity do not have such a step.
\end {proof}

\begin {remark} \label {rem-alpha-h}
  We can therefore conclude that any generalized lattice path with non-zero
  multiplicity has only two types of steps: some that go down vertically and
  others that move exactly one column to the right (with a simultaneous move up
  or down):

  \begin {center} \begin{picture}(0,0)%
\includegraphics{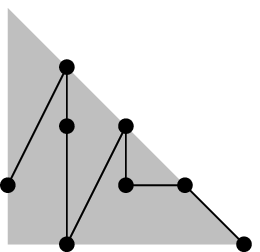}%
\end{picture}%
\setlength{\unitlength}{4144sp}%
\begingroup\makeatletter\ifx\SetFigFont\undefined%
\gdef\SetFigFont#1#2#3#4#5{%
  \reset@font\fontsize{#1}{#2pt}%
  \fontfamily{#3}\fontseries{#4}\fontshape{#5}%
  \selectfont}%
\fi\endgroup%
\begin{picture}(1150,1116)(416,-726)
\end{picture}%
 \end {center}

  For such a path we will fix the following notation for the rest of this
  section: for the vertical line $ \{ x=i \} $ we let $ h(i) $ be the
  highest $y$-coordinate of a point of $ \gamma $ on this line, and we denote
  by $ \alpha^i $ the sequence describing the lengths of the vertical steps of
  $ \gamma $ on this line. For example, for the path shown above we have
  $ h(0)=1 $, $ h(1)=3 $, $ h(2)=2 $, $ h(3)=1 $ and $ \alpha^0 = (0), \alpha^1
  =(1,1), \alpha^2=(1), \alpha^3=(0) $.
\end {remark}

\begin {proposition} \label {mult-formula}
  Let $ \alpha $ and $ \beta $ be two sequences with $ I\alpha+I\beta=d $, and
  let $\gamma$ be a generalized lattice path as above.
  \begin {enumerate}
  \item \label {mult-formula-a}
    The negative multiplicity of $\gamma$ is given by the formula
      \[ \mu_{\beta,-}(\gamma)= \sum_{(\beta^0,\ldots,\beta^d)}
         \left(\prod_{i=0}^{d-1} I^{\alpha^{i+1}+\beta^{i+1}-\beta^i} \cdot
           \binom{\alpha^{i+1}+\beta^{i+1}}{\beta^i}\right) \]
    where the sum is taken over all $(d+1)$-tuples of sequences $ (\beta^0,
    \ldots,\beta^d) $ such that $ \alpha^0+\beta^0=\beta $ and $
    I\alpha^i+I\beta^i = h(i) $ for all $i$.
  \item \label {mult-formula-b}
    The positive multiplicity of $\gamma$ is given by the formula
      \[ \mu_{\alpha,+}(\gamma)= \frac 1{I^\alpha} \cdot
         \sum_{(\beta^0,\ldots,\beta^d)}
         \left(\prod_{i=0}^{d-1} I^{\alpha^i+\beta^i-\beta^{i+1}} \cdot
           \binom{\alpha^i+\beta^i}{\beta^{i+1}}\right) \]
    where the sum is taken over all $(d+1)$-tuples of sequences $(\beta^0,
    \ldots,\beta^d)$ such that $ \beta^0=\alpha $ and $ d-i-I\beta^i =
    h(i) $ for all $i$.
  \end {enumerate}
\end {proposition}

\begin {remark} \label {rem-ilia}
  Before we give the proof of this proposition let us brief\/ly comment on
  the geometric meaning of these formulas.

  The formula of proposition \ref {mult-formula} \ref {mult-formula-a} can be
  interpreted as the number of ways to subdivide the area of $ \Delta_d $ below
  $ \gamma $ into parallelograms and triangles such that
  \begin {itemize}
  \item the subdivision contains all vertical lines $ \{ x=i \} $ below $
    \gamma $; and
  \item each triangle in the subdivision ``is pointing to the left'', i.e.\ the
    vertex opposite to its vertical edge lies to the left of this edge,
  \end {itemize}
  where each such subdivision is counted with a multiplicity equal to the
  product of the double areas of its triangles. Indeed, the sequences $ \beta^i
  $ describe the lengths of the vertical edges in the subdivisions below $
  \gamma $. The binomial coefficients $ \binom {\alpha^{i+1}+\beta^{i+1}}{
  \beta^i} $ in the formula count the number of ways to arrange the
  parallelograms and triangles, and the factors $ I^{\alpha^{i+1}+\beta^{i+1}
  -\beta^i} $ are simply the double areas of the triangles. As an example let
  us consider the path of remark \ref {rem-alpha-h} with $ \beta=(1) $. In
  this case there is only one possibility to fill the area below $ \gamma $
  with parallelograms and triangles in this way, namely as in the following
  picture on the left:

  \begin {center} \begin{picture}(0,0)%
\includegraphics{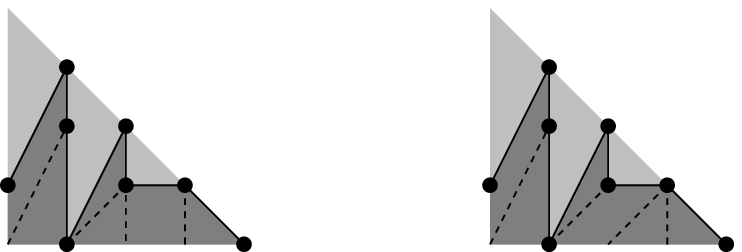}%
\end{picture}%
\setlength{\unitlength}{4144sp}%
\begingroup\makeatletter\ifx\SetFigFont\undefined%
\gdef\SetFigFont#1#2#3#4#5{%
  \reset@font\fontsize{#1}{#2pt}%
  \fontfamily{#3}\fontseries{#4}\fontshape{#5}%
  \selectfont}%
\fi\endgroup%
\begin{picture}(3355,1116)(416,-726)
\end{picture}%
 \end {center}

  (corresponding to $ \beta^0=\beta^2=\beta^3=(1), \beta^1=(0) $). As there is
  one triangle in this subdivision with double area 2 we see that $
  \mu_{\beta,-} = 2 $.

  The surprising fact about this statement is that the original definition of
  the negative multiplicity (in definition \ref {neg-mult} as well as in \cite
  {Mi03}) was set up in a way so that it also counts certain subdivisions of
  the area $ \Delta_d $ below $ \gamma $ --- but \emph {different ones}, namely
  those subdivisions that can occur as the Newton subdivisions of tropical
  curves passing through given points in a certain special position. In our
  concrete example above one can in fact see that the subdivision counted above
  \emph {does not} correspond to an actual tropical curve through the given
  points, whereas the subdivision in the picture above on the right does. One
  can therefore interpret proposition \ref {mult-formula} as the combinatorial
  statement that the number of ``column-wise'' subdivisions as described above
  agrees with the number of ``tropical'' subdivisions of the area of $ \Delta_d
  $ below $ \gamma $ (when counted with the correct multiplicities).

  It is interesting to note that for the positive multiplicity there is no such
  difference: it can be checked that the ``tropical subdivisions'' are just
  equal to the ``column-wise subdivisions'' of the area of $ \Delta_d $ above $
  \gamma $ in this case. Given this fact the statement of proposition \ref
  {mult-formula} \ref {mult-formula-b} is then almost obvious.

  In any case the nice thing about proposition \ref {mult-formula} (resp.\ the
  ``column-wise'' subdivisions) is that in these subdivisions it is easy to
  split off the first column to obtain a similar subdivision of $ \Delta_{d-1}
  $. This will be the key ingredient in the proof of the Caporaso-Harris
  formula in the lattice path set-up in theorem \ref {chlp}.
\end {remark}

\begin {proof}[Proof of proposition \ref {mult-formula}]
  We prove only part \ref {mult-formula-a} of the proposition since part \ref
  {mult-formula-b} is entirely analogous. The proof will be by induction on the
  recursive definition of $ \mu_{\beta,-} $ in definition \ref {neg-mult}. It
  is obvious that the end paths in this recursion (the paths that go from $
  (0,I\beta) $ to $ (d,0) $ along the border of $ \Delta_d $) satisfy the
  stated formula: all these paths have $ \beta^0=(0) $, so the condition $
  \alpha^0+\beta^0=\beta $ requires $ \alpha^0=\beta $, i.e.\ that the path is
  one of the paths $ \delta_\beta $ as in definition \ref {neg-mult}.

  Let us now assume that $ \gamma:[0,n] \to \Delta_d $ is an arbitrary
  generalized lattice path. By induction we know that the paths $ \gamma_-' $
  and $ \gamma_-'' $ of definition \ref {def-mu} satisfy the formula of the
  proposition. Recall that if $ k \in [1,n-1] $ is the first vertex at which $
  \gamma $ makes a right turn then $ \gamma' $ and $ \gamma'' $ are defined by
  cutting this vertex $ \gamma(k) $ (resp.\ flipping it to a parallelogram). By
  lemma \ref {no-skip} we know that $ \gamma(k-1) $ (resp.\ $ \gamma(k+1) $)
  can be at most one column to the left (resp.\ right) of $ \gamma(k) $. But
  note that $ \gamma(k-1) $ cannot be in the same column as $ \gamma(k) $ as
  otherwise the $ \lambda $-increasing path $ \gamma $ could not make a right
  turn at $ \gamma(k) $. Hence $ \gamma(k-1) $ is precisely one column left of
  $ \gamma(k) $, and we are left with two possibilities:
  \begin {itemize} \parindent 0mm \parskip 0.8ex
  \item $ \gamma(k+1) $ is in the same column $i$ as $ \gamma(k) $:

    \begin {center} \input {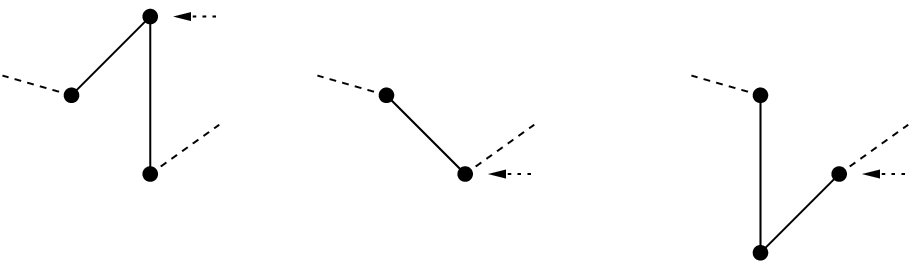} \end {center}

    Then the path $ \gamma' $ has the same values of $ h(j) $ and $ \alpha^j $
    (see remark \ref {rem-alpha-h}) as $ \gamma $, except for $ h(i) $ being
    replaced by $ h(i)-s $ and $ \alpha^i $ by $ \alpha^i-e_s $, where $s$ is
    the length of the vertical step from $ \gamma(k) $ to $ \gamma(k+1) $. So
    by induction we have
    \begin {align*}
      \mu_{\beta,-}(\gamma') &= \sum_{(\beta^0,\ldots,\beta^d)}
         \left(\prod_{j=0}^{d-1} I^{\alpha^{j+1}+\beta^{j+1}-\beta^j
           - \delta_{i,j+1} e_s} \cdot
           \binom{\alpha^{j+1}+\beta^{j+1} - \delta_{i,j+1} e_s}{\beta^j}
           \right) \\
      &= \frac 1s \cdot \sum_{(\beta^0,\ldots,\beta^d)}
         \left(\prod_{j=0}^{d-1} I^{\alpha^{j+1}+\beta^{j+1}-\beta^j} \cdot
           \binom{\alpha^{j+1}+\beta^{j+1} - \delta_{i,j+1} e_s}{\beta^j}
           \right)
    \end {align*}
    where the sum is taken over the same $ \beta $ as in the proposition. The
    path $ \gamma'' $ has the same values of $ h(j) $ and $ \alpha^j $ as $
    \gamma $ except for $ h(i) $ being replaced by $ h(i)-s $, $ \alpha^i $ by
    $ \alpha^i - e_s $, and $ \alpha^{i-1} $ by $ \alpha^{i-1}+e_s $. So by
    induction it follows that
      \[ \mu_{\beta,-}(\gamma'')= \sum_{(\beta^0,\ldots,\beta^d)}
         \left(\prod_{j=0}^{d-1} I^{\alpha^{j+1}+\beta^{j+1}-\beta^j} \cdot
           \binom{\alpha^{j+1}+\beta^{j+1} + \delta_{i-1,j+1} e_s -
           \delta_{i,j+1} e_s}{\beta^j}
           \right) \]
    where the conditions on the summation variables $ \beta^i $ are $ \alpha^0
    + \delta_{i-1,0} e_s + \beta^0 = \beta $ and $ I(\alpha^j - \delta_{i,j}
    e_s + \delta_{i-1,j} e_s) + I\beta^j = h(j)-s\delta_{i,j} $. We can make
    these conditions the same as in the proposition by replacing the summation
    variables $ \beta^{i-1} $ by $ \beta^{i-1}-e_s $, arriving at the formula
      \[ \mu_{\beta,-}(\gamma'')= \sum_{(\beta^0,\ldots,\beta^d)}
         \left(\prod_{j=0}^{d-1} I^{\alpha^{j+1}+\beta^{j+1}-\beta^j} \cdot
           \binom{\alpha^{j+1}+\beta^{j+1} - \delta_{i,j+1} e_s}{\beta^j
           - \delta_{i,j+1} e_s}
           \right). \]
    Plugging these expressions into the defining formula
      \[ \mu_{\beta,-} (\gamma) = s \cdot \mu_{\beta,-} (\gamma') +
         \mu_{\beta,-} (\gamma'') \]
    we now arrive at the desired formula of the proposition (where we simply
    use the standard binomial identity $ \binom {n-1}{k} + \binom {n-1}{k-1}
    = \binom nk $).
  \item $ \gamma(k+1) $ is one column right of $ \gamma(k) $: The idea here is
    the same as in the previous case. The actual computation is simpler
    however as the path $ \gamma' $ does not give a contribution by lemma \ref
    {no-skip}.

    \begin {center} \input {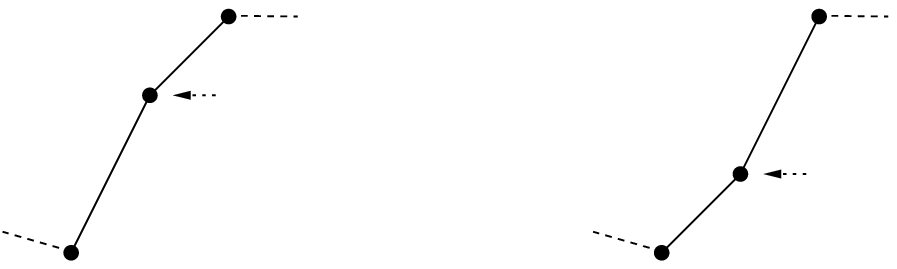} \end {center}

    The path $ \gamma'' $ has the same values of $ h(j) $ and $ \alpha^j $ as $
    \gamma $, except for $ h(i) $ being replaced by $ h(i)-s+t $, where $i$ is
    the column of $ \gamma(k) $, and $s$ and $t$ are the vertical lengths of
    the steps before and after $ \gamma(k) $. So by induction we have
      \[ \mu_{\beta,-}(\gamma'')= \sum_{(\beta^0,\ldots,\beta^d)}
         \left(\prod_{j=0}^{d-1} I^{\alpha^{j+1}+\beta^{j+1}-\beta^j} \cdot
           \binom{\alpha^{j+1}+\beta^{j+1}}{\beta^j}
           \right) \]
    where the conditions on the $ \beta^j $ are $ \alpha^0+\beta^0=\beta $ and
    $ I\alpha^j + I\beta^j = h(j) - (s-t) \delta_{i,j} $. Note that $ s-t>0 $
    since $ \gamma $ makes a right turn. As in the previous case we can make
    the conditions on the $ \beta^j $ the same as in the proposition by
    replacing the summation variables $ \beta^i $ by $ \beta^{i-1}+
    \alpha^{i+1}+\beta^{i+1}-\beta^i $. We then arrive at
    \begin{align*}
      \mu_{\beta,-}(\gamma'')
        =& \sum_{(\beta^0,\ldots,\beta^d)} 
           \prod_{j=0}^{d-1} \Bigg( I^{\alpha^{j+1}+\beta^{j+1}-\beta^j} \cdot
           \\
         & \binom{\alpha^{j+1}+\beta^{j+1}+ \delta_{i,j+1}
           (\beta^{i-1}+\alpha^{i+1}+\beta^{i+1}-2\beta^i)}{\beta^j+
           \delta_{i,j} (\beta^{i-1}+\alpha^{i+1}+\beta^{i+1}-2\beta^i)}
           \Bigg).
    \end {align*}
    This is already the formula of the proposition except for the factors
      \[ \binom {\alpha^i+\beta^i}{\beta^{i-1}}
         \binom {\alpha^{i+1}+\beta^{i+1}}{\beta^i} \]
    being replaced by
    \begin{align*}
        & \binom {\alpha^i+\beta^i+(\beta^{i-1}+\alpha^{i+1}+\beta^{i+1}
          -2\beta^i)}{\beta^{i-1}} \\
        & \qquad \cdot \binom {\alpha^{i+1}+\beta^{i+1}}{\beta^i +
          (\beta^{i-1}+\alpha^{i+1}+\beta^{i+1}-2\beta^i)}.
    \end{align*}
    But these terms are the same by the identity $ \binom
    {n+k+l}{n+k} \binom {n+k}{n} = \binom {n+k+l}{n+l} \binom {n+l}{n}
    $ with $ n=\beta^{i-1} $, $ k=\beta^i-\beta^{i-1} $, and $ l=\alpha^{i+1}+
    \beta^{i+1}-\beta^i$ (note that $ \alpha^i = (0) $ for our path $ \gamma
    $).
  \end {itemize}
  This proves the proposition.
\end {proof}

\begin{remark} \label{rem-gen}
  Note that it is important for the second step of the proof above that the two
  boundary lines of $\Delta_d$ below and above $\gamma$ --- the line $\{y=0\}$
  respectively the diagonal line from $(0,d)$ to $(d,0)$ --- are indeed
  straight lines. We cannot generalize the proof to polygons which contain a
  vertex above respectively below $\gamma$, because then the heights of the
  three columns of $\gamma(k-1)$, $\gamma(k)$ and $\gamma(k+1)$ cannot be
  described as $h(i)$, $h(i)+s$ and $h(i)+s+t$. So we cannot use the identity $
  \binom {n+k+l}{n+k} \binom {n+k}{n} = \binom {n+k+l}{n+l} \binom {n+l}{n} $.
  The picture below shows a polygon for which the formula of proposition
  \ref{mult-formula} does not hold. The column-wise subdivisions would predict
  $0$ as negative multiplicity for the path. However, the path $\gamma''$ is a
  valid end path, so we get $1$.

  \begin {center} \begin{picture}(0,0)%
\includegraphics{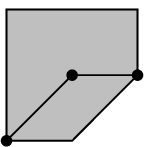}%
\end{picture}%
\setlength{\unitlength}{4144sp}%
\begingroup\makeatletter\ifx\SetFigFont\undefined%
\gdef\SetFigFont#1#2#3#4#5{%
  \reset@font\fontsize{#1}{#2pt}%
  \fontfamily{#3}\fontseries{#4}\fontshape{#5}%
  \selectfont}%
\fi\endgroup%
\begin{picture}(659,641)(766,-330)
\end{picture}%
 \end {center}

  The proof can be generalized to polygons where the boundaries above and below
  $\gamma$ are straight lines, for example to rectangles (see \ref{ch-rect}).
\end{remark}

We are now ready to prove the Caporaso-Harris formula in the lattice path
set-up:

\begin{theorem} \label{chlp}
  The numbers $ N_{\path}^{\alpha,\beta}(d,g) $ satisfy the Caporaso-Harris
  formula. In particular we have $ N_{\path}^{\alpha,\beta}(d,g) =
  N_{\cplx}^{\alpha,\beta}(d,g) $ for all $ d,g,\alpha,\beta $.
\end{theorem}

\begin {proof}
  The idea of the proof is simply to list the possibilities of the first step
  of the path $\gamma$. Let $\gamma$ be a $\lambda$-increasing path from
  $(0,I\beta)$ to $q=(d,0)$. Then we have one of the following two cases:

  Case 1: The point $\gamma(1)$ is on the line $ \{x=0\} $. Then
  $ \gamma(1) $ must be $(0,I\beta-k)$ for some $\beta_k \neq 0$ as
  otherwise $\mu_{\beta,-}(\gamma)$ would be $0$. It follows that $ \gamma
  \mid_{[1,2d+g+|\beta|-1]}$ is a path from $(0,d-I(\alpha+e_k))$ and
  with $ \mu_{\alpha,\beta}(\gamma)= k\cdot\mu_{\alpha+e_k,\beta-e_k}
  (\gamma\mid_{[1,2d+g+|\beta|-1]})$. Therefore the paths $ \gamma $ with
  $ \gamma(1) \in s $ contribute
    \[ \sum_{k:\beta_k>0} k\cdot N_{\path}^{\alpha+e_k,\beta-e_k}(d,g)
       \]
  to the number $ N_{\path}^{\alpha,\beta}(d,g) $.

  Case 2: The point $\gamma(1)$ is not on $\{x=0\}$. Then it must be on the
  line $\{x=1\}$ by lemma \ref {no-skip}. From proposition \ref {mult-formula}
  it follows that both the negative and the positive multiplicity can be
  computed as a product of a factor coming from the first column and the
  (negative resp.\ positive) multiplicity of the restricted path $ \tilde
  \gamma := \gamma|_{[1,2d+g+|\beta|-1]} $. More precisely, we have
  \begin {align*}
    \mu_{\alpha,\beta} (\gamma)
    &= \mu_{\beta,-} (\gamma) \cdot \mu_{\alpha,+} (\gamma) \\
    &= \sum_{\beta'} I^{\beta'-\beta} \binom {\beta'}{\beta} \mu_{\beta',-}
       (\tilde \gamma) \cdot \sum_{\alpha'}
       \binom {\alpha}{\alpha'} \mu_{\alpha',+}(\tilde \gamma) \\
    &= \sum_{\alpha',\beta'} I^{\beta'-\beta} \binom {\beta'}{\beta}
       \binom {\alpha}{\alpha'} \cdot \mu_{\alpha',\beta'} (\tilde \gamma).
  \end {align*}
  So the contribution of the paths with $ \gamma(1) \notin s $ to $
  N_{\path}^{\alpha,\beta}(d,g) $ is
    \[ \sum I^{\beta'-\beta} \binom {\beta'}{\beta} \binom {\alpha}{\alpha'}
         \cdot N_{\path}^{\alpha',\beta'}(d-1,g') \]
  where the sum is taken over all possible $\alpha',\beta'$ and $g'$. Let us
  figure out what these possible values are. It is clear that $ \alpha' \le
  \alpha $ and $ \beta \le \beta' $. Also, $I\alpha'+I\beta'=d-1$ must be
  fulfilled. As $ \tilde \gamma $ has one step less than $\gamma$ we know that
  $2d+g+|\beta|-1-1=2(d-1)+g'+|\beta'|-1 $ and hence $ g-g'=|\beta'-\beta|-1$.
  A path $\epsilon:[0,n]\rightarrow \Delta$ from $(0,I\beta)$ to $q$ that meets
  all lattice points of $\Delta$ has $|\beta|+d (d+1)/2$ steps. As $\gamma$ has
  $2d+g-1+|\beta|$ steps, $|\beta|+d(d+1)/2 -(2d+g-1+|\beta|)=(d-1)(d-2)/2-g$
  lattice points are missed by $\gamma$. But $ \tilde \gamma $ cannot miss more
  points, therefore $(d-2)(d-3)/2-g' \leq (d-1)(d-2)/2-g $, i.e.\ $ d-2 \geq
  g-g'$.
\end {proof}
  
\begin{remark} \label{ch-rect}
  The same argument can also be applied to other polygons $ \Delta $. For
  example, the analogous recursion formula for $\P^1 \times \P^1$, i.e.\ for a
  $ d' \times d $ rectangle $\Delta_{(d',d)}$ reads
  \begin{align*}
    N_{\path}^{\alpha,\beta}((d',d),g)
      =& \sum_{k:\beta_k>0} k\cdot N_{\path}^{\alpha+e_k,\beta-e_k}
        ((d',d),g) \\
      &+ \sum I^{\beta'-\beta} \binom{\alpha}{\alpha'} \binom{\beta'}{\beta}
        N_{\path}^{\alpha',\beta'}((d'-1,d),g')
  \end{align*}
  for all $ \alpha ,\beta $ with $ I\alpha+I\beta=d $, where the second sum is
  taken over all $\alpha',\beta',g'$ such that $\alpha \leq \alpha$, $\beta'
  \geq \beta$, $I\alpha'+I\beta'=d$, $g-g'\leq d-1$ and $|\beta-\beta'|=g'-g-1
  $.
\end{remark}


\section {The Caporaso-Harris formula in the tropical set-up}
  \label{sec-tropical}

\subsection {Relative Gromov-Witten invariants in tropical geometry}

We will now define the analogues of the numbers $ N_{\cplx}^{\alpha,\beta}(d,g)
$ in terms of tropical curves. The definition is quite straightforward:

\begin{definition} \label{ntrop}
  Let $C$ be a simple tropical curve of degree $\Delta_d$ and genus $g$ with
  $\alpha_i$ fixed and $\beta_i$ non-fixed unbounded ends to the left of weight
  $i$ for all $i$. We define the $(\alpha,\beta)$-\df{multiplicity} of $C$ to
  be
    \[ \mult_{\alpha,\beta}(C):=\frac{1}{I^{\alpha}} \cdot \mult(C) \]
  where $ \mult(C) $ is the usual multiplicity as in definition \ref{tropdef}
  \ref{tropdef-f}.

  Let $ d \ge 0 $ and $g$ be integers, and let $ \alpha $ and $ \beta $ be
  sequences with $ I\alpha+I\beta=d $. Then we define $
  N_{\trop}^{\alpha,\beta}(d,g)$ to be the number of tropical curves of
  degree $\Delta_d$ and genus $g$ with $\alpha_i$ fixed and $\beta_i$ non-fixed
  unbounded ends to the left of weight $i$ for all $i$ that pass in addition
  through a set $\mathcal{P}$ of $ 2d+g+|\beta|-1 $ points in general position.
  The curves are to be counted with their respective $(\alpha,\beta)
  $-multiplicities. By \cite{GM051} this definition does not depend on the
  choice of marked points and fixed unbounded ends.
\end{definition}

First of all we will show that this definition actually coincides with the
lattice path construction of definition \ref{nd}:

\begin{theorem} \label{npantrop}
  For all $ d,g,\alpha,\beta $ we have $N_{\trop}^{\alpha,\beta}(d,g)=
  N_{\path}^{\alpha,\beta}(d,g)$.
\end{theorem}

\begin{proof}
  The proof is analogous to the proof of \cite{Mi03} theorem 2. As usual we
  choose $\lambda(x,y)=x-\varepsilon y$. Let $\mathcal{P}$ be a set of
  $2d+g+|\beta|-1$ points on a line $H$ orthogonal to the kernel of $\lambda$
  such that the distance between $p_i$ and $p_{i+1}$ is much bigger than the
  distance of $p_{i-1}$ and $p_i$ for all $i$, and such that all points lie
  below the fixed ends. In other words, if the fixed ends have the
  $y$-coordinates $y_{1},\ldots,y_{|\alpha|} $ then the $y$-coordinates of
  $p_i$ are chosen to be less than all $y_{1},\ldots,y_{|\alpha|}$. Our aim is
  to show that the number of tropical curves through this special configuration
  is equal to the number of lattice paths as in section \ref{sec-lattice}. Let
  $C$ be a tropical curve with the right properties through this set of points.
  Mark the points where $H$ intersects the fixed ends. The proof of
  \cite{Mi03} theorem 2 tells us that the edges of $ \Delta $ (dual to the
  edges of the curves where they meet $\mathcal{P}$ and the new marked points)
  form a $\lambda$-increasing path from $p=(0,d)$ to $q=(d,0)$. The fact that
  the fixed ends lie above all other points tells us that the path starts with
  $\alpha_i$ steps of lengths $i$. So we can cut the first part and get a path
  from $(0,I\beta)$ to $q$ with the right properties. 

  Next, let a path $\gamma:[0,2d+g+|\beta|-1]\rightarrow \Delta_d$ be given
  that starts at $ (0,d-I\alpha) $ and ends at $q$. Extend $\gamma$ to a path
  $\gamma_{\alpha}:[0,|\alpha|+n]\rightarrow \Delta_d$ by adding $\alpha_i$
  steps of integer length $i$ at $\{x=0\}$ from $\gamma_{\alpha}(0)=p$ to
  $\gamma_{\alpha}(|\alpha|)=(0,I\beta)$. Add the steps of integer lengths $i$
  in an order corresponding to the order of the fixed ends. The recursive
  definition of $\mu_{\beta,-}(\gamma_{\alpha})$ corresponds to counting the
  possibilities for a dual tropical curve in the half plane below $H$ through
  $\mathcal{P}$. Passing from $\gamma_{\alpha}$ to $\gamma_{\alpha}'$ and
  $\gamma_{\alpha}''$ corresponds to counting the possibilities in a strip
  between $H$ and a parallel shift of $H$. We end up with a path $\delta_-$
  which begins with $\alpha_i$ steps of length $i$ and continues with $\beta_i$
  steps of lengths $i$. This shows that all possible dual curves have the right
  horizontal ends. Furthermore, $\mu_{\beta,-}(\gamma_{\alpha})$ coincides with
  the number of possible combinatorial types of the curve in the half plane
  below $H$ times the multiplicity of the part of the curve in the half plane
  below $H$. With the same arguments we get that $\mu_{\alpha,+}(\gamma_{
  \alpha})$ is equal to the number of possibilities for the combinatorial type
  times the multiplicity in the upper half plane. Altogether, we have
  \begin{align*}
    N_{\path}^{\alpha,\beta}(d,g)
      &= \sum_\gamma \mult_{\alpha,\beta}(\gamma) \\
      &=  \frac{1}{I^{\alpha}} \sum_\gamma \mult_{\beta,-}(\gamma_{\alpha})
        \cdot \mult_{\alpha,+}(\gamma_{\alpha})\\
      &= \frac{1}{I^{\alpha}}\sum_C \mult(C)=\sum_C \mult_{\alpha,\beta}(C)\\
      &= N_{\trop}^{\alpha,\beta}(d,g),
  \end{align*}
  where $C$ runs over all tropical curves with the right properties and
  $\gamma$ runs over all paths with the right properties.
\end{proof}

\subsection{The Caporaso-Harris formula} \label {ch-red}

Of course it follows from theorems \ref{chlp} and \ref{npantrop} that the
numbers $ N_{\trop}^{\alpha,\beta} (d,g) $ satisfy the Caporaso-Harris formula.
In this section we will prove this statement directly without using Newton
polyhedra and lattice paths. The advantage of this method is that it may be
possible to generalize it to curves in higher-dimensional varieties (where the
concept of Newton polyhedra cannot be used to describe tropical curves).

\begin{theorem} \label{trCH}
  The numbers $ N_{\trop}^{\alpha,\beta}(d,g)$ satisfy the Caporaso-Harris
  formula.
\end{theorem}

\begin{proof}
  Let $ \varepsilon>0 $ be a small and $ N>0 $ a large real number. We
  choose the fixed unbounded left ends and the set $ \mathcal {P} =
  \{p_1,\dots,p_n \} $ in tropical general position so that
  \begin {itemize}
  \item the $y$-coordinates of all $ p_i $ and the fixed ends are in the open
    interval $ (-\varepsilon,\varepsilon) $;
  \item the $x$-coordinates of $ p_2,\dots,p_n $ are in the open interval
    $ (-\varepsilon,\varepsilon) $;
  \item the $x$-coordinate of $p_1$ is less than $ -N $.
  \end {itemize}
  In other words, we keep all conditions for the curves in a small horizontal
  strip and move $ p_1 $ very far to the left.

  \begin {center} \input {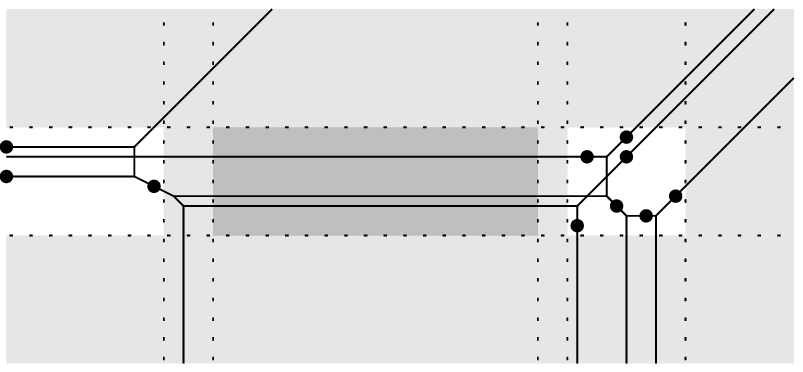} \end {center}

  Let us consider a tropical curve $C$ satisfying the given conditions. We want
  to show that $C$ must ``look as in the picture above'', i.e.\ that the curve
  ``splits'' into two parts joined by only horizontal lines.

  First of all we claim that no vertex of $C$ can have its $y$-coordinate below
  $ -\varepsilon $: otherwise let $V$ be a vertex with lowest $y$-coordinate.
  By the balancing condition there must be an edge pointing downwards from $V$.
  As there is no vertex below $V$ this must be an unbounded edge. For degree
  reasons this edge must have direction $ (0,-1) $ and weight 1, and it must be
  the only edge pointing downwards. By the balancing condition it then follows
  that at least one other edge starting at $V$ must be horizontal. Again by the
  balancing condition we can then go from $V$ along this horizontal edge to
  infinity in the region $ \{ y \le -\varepsilon \} $. As there are no marked
  points in this region we could go from $V$ to infinity without passing
  marked points in two different ways, which is a contradiction to remark
  \ref{rem-strings}. (This is an important argument which we will use several
  times in this proof.)

  In the same way we see that no vertex of $C$ can have its $y$-coordinate
  above $ \varepsilon $.

  Next, consider the rectangle
    \[ R := \{(x,y);\; -N \le x \le -\varepsilon,
         -\varepsilon \le y \le \varepsilon \}. \]
  We want to study whether there can be vertices of $C$ within $R$. Let $ C_0 $
  be a irreducible component of $ C \cap R $. Note that any end of $ C_0 $
  leaving $R$ at the top or bottom edge must go straight to infinity as we have
  just seen that there are no vertices of $C$ above or below $R$. If there are
  ends of $ C_0 $ leaving $R$ at the top \emph{and} at the bottom then we could
  again go from one infinite end of $C$ to another without passing a marked
  point, again in contradiction to remark \ref{rem-strings}. So we may assume
  without loss of generality that $ C_0 $ does not meet the top edge of $R$.
  With the same argument, we can see that $C_0$ can meet the top edge of $R$
  only in one point.

  It follows by the balancing condition that all edges of $ C_0 \cap R $ that
  are not horizontal project to the $x$-axis to a union of two (maybe empty)
  intervals $ [-N,x_1] \cup [x_2,-\varepsilon] $. But note that the number of
  edges of $C$ as well as the minimum slope of an edge (and hence the maximum
  distance an edge can have within $R$) are bounded by a constant that depends
  only on the degree of the curves. So we can find $ a,b \in \R $ (that depend
  only on $d$) such that the interval $ [a,b] $ is disjoint from $ [-N,x_1]
  \cup [x_2,-\varepsilon] $, or in other words such that there are no
  non-horizontal edges in $ [a,b] \times [-\varepsilon,\varepsilon] $. In
  particular, there are then no vertices of $C$ in $ [a,b] \times \R $. Hence
  we see that the curve $C$ must look as in the picture above: we can
  ``cut'' it at any line $ x=x_0 $ with $ a < x_0 < b $ and obtain curves on
  both sides of this line that are joined only by horizontal lines.

  There are now two cases to distinguish:
  \begin {enumerate} \parindent0mm \parskip 0.5ex
  \item $p_1$ lies on a horizontal non-fixed end of $C$. Then the region
    where $ x \le -\varepsilon $ consists of only horizontal lines. (Otherwise
    we could again go from one unfixed end to another without meeting a marked
    point, in contradiction to remark \ref{rem-strings}.) We can hence consider
    $C$ as having one more fixed end at $ p_1 $ and passing through
    $\mathcal{P} \setminus \{p_1\}$. We just have to multiply with the weight
    of this end, as the multiplicity of curves with fixed ends is defined as
    $\frac{1}{I^{\alpha}}\cdot \mult(C)$. Therefore the contribution of these
    curves to $ N_{\trop}^{\alpha,\beta}(d,g) $ is
      \[ \sum_{k:\beta_k>0} k\cdot N_{\trop}^{\alpha+e_k,\beta-e_k}(d,g). \]
  \item $p_1$ does not lie on a horizontal end of $C$ (as in the picture
    above). Then $C$ can be separated into a left and a right part as above.
    As the left part contains only one marked point it follows again by remark
    \ref{rem-strings} that the left part has exactly one end in direction
    $ (0,-1) $ and $ (1,1) $ each, together with some more horizontal ends.
    
    Hence the curve on the right must have degree $ d-1 $. Let us denote this
    curve by $C'$.

    How many possibilities are there for $C$? Assume that $\alpha' \leq
    \alpha$ of the fixed horizontal ends only intersect the part $C \setminus
    C'$ and are not adjacent to a $3$-valent vertex of $C\setminus C'$. Then
    $C'$ has $\alpha'$ fixed horizontal ends. Given a curve $C'$ of degree
    $d-1$ with $\alpha'$ fixed ends through $\mathcal{P}\setminus \{p_1\}$,
    there are $\binom{\alpha}{\alpha'}$ possibilities to choose which fixed
    ends of $C$ belong to $C'$. $C'$ has $d-1-I\alpha'$ non-fixed
    horizontal ends. Let $\beta'$ be a sequence which fulfills $I\beta'=d-1-
    I\alpha'$, hence a possible choice of weights for the non-fixed ends of
    $C'$. Assume that $\beta''\leq \beta'$ of these ends are adjacent to a
    $3$-valent vertex of $C\setminus C'$ whereas $\beta'-\beta''$ ends
    intersect $C\setminus C'$. The irreducible component of $C\setminus C'$
    which contains $p_1$ has to contain the two ends of direction $(0,-1)$ and
    $(1,1)$ due to the balancing condition. Also, it can contain some ends of
    direction $(-1,0)$ --- but these have to be fixed ends then, as $p_1$
    cannot separate more than two (nonfixed) ends. So all the $\beta$ nonfixed
    ends of direction $(-1,0)$ have to intersect $C\setminus C'$ --- therefore
    they have to be ends of $C'$. That is, $\beta'-\beta''= \beta$ (in
    particular $\beta' \geq \beta$). Given $C'$, there are $ \binom{\beta'}{
    \beta}$ possibilities to choose which ends of $C'$ are also ends of $C$.
    Furthermore, we have
    \begin {align*}
      \mult_{\alpha,\beta}(C)
        &= \frac{1}{I^{\alpha}} \mult(C)
         = \frac{1}{I^{\alpha}} \cdot I^{\alpha-\alpha'}\cdot
           I^{\beta'-\beta}\cdot \mult(C') \\
        &= \frac{1}{I^{\alpha'}}\cdot I^{\beta'-\beta}\cdot \mult(C')
         = I^{\beta'-\beta}\cdot \mult_{\alpha',\beta'}(C')
    \end{align*}
    where the factors $I^{\alpha-\alpha'}$ and $I^{\beta'-\beta}$ arise due to
    the $3$-valent vertices which are not part of $C'$.

    To determine the genus $g'$ of $C'$, note that $C'$ has by $
    |\alpha+\beta''|$ fewer vertices than $C$ and by $|\alpha+\beta''|-1+
    |\beta''|$ fewer bounded edges --- there are $|\alpha+\beta''|-1$ bounded
    edges in $C\setminus C'$, and $ |\beta''|$ bounded edges are cut. Hence
    $g'=1-\#\Gamma^0+|\alpha+\beta''|-\#\Gamma^1_0-|\alpha+\beta''|-|\beta''|=
    g-(|\beta''|-1)$. Furthermore, $g-g' \leq d-2$ as at most $d-2$ loops may
    be cut. Now given a curve $C'$ with $\alpha'$ fixed and $\beta'$ nonfixed
    bounded edges through $\mathcal{P}\setminus \{p_1\}$, and having chosen
    which of the $\alpha$ fixed ends of $C$ are also fixed ends of $C'$ and
    which of the $\beta'$ ends of $C'$ are also ends of $C$, there is only one
    possibility to add an irreducible component through $p_1$ to make it a
    possible curve $C$ with $\alpha$ fixed ends and $\beta$ nonfixed. The
    positions and directions of all bounded edges are prescribed by the point
    $p_1$, by the positions of the $\beta'-\beta$ ends to the left of $C'$, and
    by the $\alpha-\alpha'$ fixed ends. Hence we can count the possibilities
    for $C'$ (times the factor $\binom{\alpha}{\alpha'}\cdot \binom{\beta'}{
    \beta}\cdot I^{\beta'-\beta}$) instead of the possibilities for $C$ (where
    the possible choices for $\alpha'$, $\beta'$ and $g'$ have to satisfy just
    the conditions we know from the Caporaso-Harris formula). 
  \end {enumerate}
  The sum of these two contributions gives the required Caporaso-Harris
  formula.
\end{proof}

\subsection{The tropical Caporaso-Harris formula for irreducible curves}
  \label{ch-irr}

So far we have only considered not necessarily irreducible curves since this
case is much simpler (the irreducibility condition is hard to keep track of
when we split a curve into two parts as in the proof of theorem \ref{trCH}).
In this section we want to show how the ideas of section \ref{ch-red} can be
carried over to the case of irreducible curves.

\begin{definition} \label{ntropirr}
  Let $N_{\trop}^{\irr,(\alpha,\beta)}(d,g)$ be the number of
  irreducible tropical curves of degree $\Delta_d$ and genus $g$ with
  $\alpha_i$ fixed and $\beta_i$ non-fixed horizontal left ends of weight $i$
  for all $i$ that pass in addition through a set of $|\beta|+2d+g-1$ points
  in general position. Again, the curves are to be counted with their
  $(\alpha,\beta)$-multiplicity as of definition \ref{ntrop}. As in definition
  \ref{ntrop} it follows from \cite{GM051} that these numbers do not depend on
  the choice of the points and the positions of the fixed ends.
\end{definition}

\begin{theorem}
  The numbers $N_{\trop}^{\irr,(\alpha,\beta)}(d,g)$ satisfy the recursion
  relations
  \begin{align*}
    N_{\trop}^{\irr,(\alpha,\beta)}(d,g)
      =& \sum_{k:\beta_k>0} k\cdot N_{\trop}^{\irr,(\alpha+e_k,\beta-e_k)}
        (d,g) \\
      & \quad + \sum \frac{1}{\sigma} \binom{2d+g+|\beta|-2}{2d_1+g_1+|\beta^1|
        -1,\ldots,2d_k+g_k+|\beta^k|-1} \\
      & \qquad \cdot \binom{\alpha}{\alpha^1,\ldots,\alpha^k}  \\
      &\qquad \cdot \prod_{j=1}^k \left(\binom {\beta^j}{\beta^j-\beta^{j'}}
        \cdot I^{\beta^{j'}} \cdot N_{\trop}^{\irr,(\alpha^j,\beta^j)}(d_j,
        g_j)\right)
  \end{align*}
  with the second sum taken over all collections of integers $d_1,\ldots,
  d_k$ and $g_1,\ldots,g_k$ and all collections of sequences $\alpha^1,\ldots,
  \alpha^k$, $\beta^1,\ldots,\beta^k$ and $\beta^{1'},\ldots,\beta^{k'}$
  satisfying 
  \begin{align*}
    \alpha^1+\ldots+\alpha^k &\leq \alpha; \\
    \beta^1+\ldots+\beta^k &= \beta+\beta^{1'}+\ldots+\beta^{k'}; \\
    |\beta^{j'}| &> 0; \\
    d_1+\ldots+d_k &= d-1; \\
    g-(g_1+\ldots+g_k) &= |\beta^{1'}+\ldots+\beta^{k'}|+k.
  \end{align*}
  Here as usual $\binom{n}{a_1,\ldots,a_k}$ denotes the multinomial coefficient
  \begin{displaymath}
    \binom{n}{a_1,\ldots,a_k}=\frac{n!}{a_1!\cdot\ldots \cdot
      a_k!(n-a_1-\ldots-a_k)!}
  \end{displaymath}
  and correspondingly, for sequences $\alpha,\alpha^1,\ldots,\alpha^k$ the
  multinomial coefficient is
  \begin{displaymath}
    \binom{\alpha}{\alpha^1,\ldots,\alpha^k}=\prod_i
      \binom{\alpha_i}{\alpha^1_i,\ldots,\alpha^k_i}.
  \end{displaymath} 
  The number $\sigma$ is defined as follows: Define an equivalence relation on
  the set $\{1,2,\ldots,k\}$ by $i: \sim j$ if $d_i=d_j$, $g_i=g_j$, $\alpha^i=
  \alpha^j$, $\beta^i=\beta^j$ and $\beta^{i'}=\beta^{j'}$. Then $\sigma$ is
  the product of the factorials of the cardinalities of the equivalence
  classes. 
\end{theorem}

Note that this recursion formula coincides with the corresponding
Caporaso-Harris formula for irreducible curves (see \cite{CH98} section 1.4).

\begin{proof}
  The proof is analogous to that of theorem \ref{trCH}. Fix a set $\mathcal{P}
  =\{p_1,\ldots,p_n\}$ in tropical general position with $ p_1 $ very far left
  as in the proof of theorem \ref{trCH}. Let $C$ be an irreducible tropical
  curve with the right properties through the points. The first term in the
  recursion formula (that corresponds to curves with only horizontal lines in
  the area where $ x \le -\varepsilon $) follows in the same way as in theorem
  \ref{trCH}. So assume that $ p_1 $ does not lie on a horizontal end of $C$.
  As before we get a curve $C'$ of degree $ d-1 $ to the right of the cut. The
  curve $C'$ does not need to be irreducible however. It can split in $k$
  irreducible components $C_1,\ldots,C_k$ of degree $d_1,\ldots,d_k$. Then
  $d_1+\ldots+d_k=d-1$. As before, we would like to count the possibilities
  for the $C_j$ separately, and then determine how many ways there are to make
  a possible curve $C$ out of a given choice of $C_1,\ldots, C_k$. So let the
  $C_j$ be curves of degree $d_j$ through the set of points $\mathcal{P}
  \setminus \{p_1\}$. Let $C_j$ have $\alpha^j$ fixed horizontal ends and
  $\beta^j$ nonfixed horizontal ends, satisfying $|\beta^j|+|\alpha^j|=d_j$.
  Then $\alpha^1+\ldots+\alpha^k \leq \alpha$, and there are $\binom{\alpha}{
  \alpha^1,\ldots,\alpha^k}$ possibilities to choose which fixed ends of $C$
  belong to which $C_j$. As before, the irreducible component of $C\setminus
  (C_1\cup \ldots \cup C_k)$ which contains $p_1$ is fixed by only one point,
  therefore it contains none of the $\beta$ nonfixed ends of $C$. That is, all
  $\beta$ nonfixed ends have to intersect the part $C\setminus (C_1\cup \ldots
  \cup C_k)$ and have to be ends of one of the $C_j$ then. Assume that $
  \beta^{j'}$ of the $\beta^j$ ends are adjacent to a $3$-valent vertex of
  $C\setminus (C_1\cup \ldots \cup C_k)$ whereas $\beta^j-\beta^{j'}$ ends
  just intersect $C\setminus (C_1\cup \ldots \cup C_k)$. As $C$ is irreducible
  we must have $|\beta^{j'}|>0$ as otherwise $C_j$ would form a separate
  component of $C$. Then given the curves $C_j$ through $\mathcal{P}
  \setminus \{p_1\}$, there are $\binom{\beta^j}{\beta^j-\beta^{j'}}$
  possibilities to choose which of the nonfixed ends of $C_j$ are also ends of
  $C$, and $\beta+\sum \beta^{j'}=\sum \beta^j$. Each $C_j$ is fixed by
  $2d_j+g_j+ |\beta^j|-1$ points, where $g_j$ denotes the genus of $C_j$.
  (There cannot be fewer points on one of the $C_j$, since otherwise the
  unbounded ends or loops could not be separated by the points, in
  contradiction to remark \ref{rem-strings}.) Therefore, there are
  \begin{displaymath}
    \binom{2d+g+|\beta|-2}{2d_1+g_1+|\beta^1|-1,\ldots,2d_k+g_k+|\beta^k|-1}
  \end{displaymath}
  possibilities to distribute the points $p_2,\ldots,p_n$ on the components
  $C_1,\ldots, C_k$. Furthermore, we have
  \begin{align*}
    \mult_{\alpha,\beta}(C)
      &= \frac{1}{I^{\alpha}} \mult(C) \\
      &= \frac{1}{I^{\alpha}} \cdot I^{\alpha-\alpha^1-\ldots-\alpha^k}
           \cdot I^{\beta^{1'}+\ldots+\beta^{k'}}
           \cdot \mult(C_1)\cdot \ldots \cdot \mult(C_k) \\
      &= \frac{1}{I^{\alpha^1+\ldots+\alpha^k}}
           \cdot I^{\beta^{1'}+\ldots+\beta^{k'}}
           \cdot \mult(C_1)\cdot \ldots \cdot \mult(C_k) \\
      &= I^{\beta^{1'}+\ldots+\beta^{k'}}
           \cdot \mult_{\alpha^1,\beta^1}(C_1)\cdot \ldots
           \cdot \mult_{\alpha^k,\beta^k}(C_k)
  \end{align*}
  where the factors $I^{\alpha-\alpha^1-\ldots-\alpha^k}$ and $I^{\beta^{1'}+
  \ldots+\beta^{k'}}$ arise due to the $3$-valent vertices which are not part
  of $C_1,\ldots,C_k$. Concerning the genus, note that $C$ has $|\alpha-
  \alpha^1-\ldots-\alpha^k+\beta^{1'}+\ldots+\beta^{k'}|$ more vertices than
  $C_1\cup \ldots \cup C_k$ and $|\alpha-\alpha^1-\ldots-\alpha^k+\beta^{1'}+
  \ldots+\beta^{k'}|-1+|\beta^{1'}+\ldots+\beta^{k'}|$ more bounded edges
  (there are $|\alpha-\alpha^1-\ldots-\alpha^k+\beta^{1'}+\ldots+\beta^{k'}|-1$
  bounded edges in $C\setminus (C_1\cup \ldots \cup C_k)$ and $|\beta^{1'}+
  \ldots+\beta^{k'}|$ bounded edges are cut). Hence
  \begin{align*}
    g &= 1+ g_1+\ldots+g_k -k \\
      &\qquad -(|\alpha-\alpha^1-\ldots-\alpha^k+\beta^{1'}+\ldots+\beta^{k'}|) \\
      &\qquad + |\alpha-\alpha^1-\ldots-\alpha^k+\beta^{1'}+\ldots+\beta^{k'}|
          -1+|\beta^{1'}+\ldots+\beta^{k'}| \\
      &= \sum g_j + \sum |\beta^{j'}|-k.
  \end{align*}
  This proves the recursion formula except for the factor $\frac{1}{\sigma}$.
  This factor is simply needed because up to now we count different curves if
  two components $C_i$ and $C_j$ of $C'$ are identical, depending on whether
  $C_i$ is the $i$-th component or $C_j$ is the $i$-th component. Therefore we
  have to divide by $\sigma$.
\end{proof}

\providecommand{\bysame}{\leavevmode\hbox to3em{\hrulefill}\thinspace}
\providecommand{\MR}{\relax\ifhmode\unskip\space\fi MR }
\providecommand{\MRhref}[2]{%
  \href{http://www.ams.org/mathscinet-getitem?mr=#1}{#2}
}
\providecommand{\href}[2]{#2}

\end{document}